\documentclass[fleqn]{article}
\usepackage{latexsym}
\usepackage{amssymb}
\usepackage{amsmath}
\usepackage{color}
\usepackage{graphicx}
\usepackage{mathrsfs}
\allowdisplaybreaks

\newtheorem{theorem}{\color{black} Theorem}%[section]
\newtheorem{lemma}{\color{black} Lemma}[section]
\newtheorem{proposition}{\color{black} Proposition}[section]
\newtheorem{definition}{\color{black} Definition}[section]
\newtheorem{remark}{\color{black} Remark}[section]

\newtheorem{example}{\color{black} Example}[section]

\begin{document}

%\large
\title{\bf The random case of Conley's theorem\thanks{Published in:
{\it Nonlinearity} {\bf 19} (2006), 277-291.}}

\author{
{\bf Zhenxin Liu}\thanks{{\it E-mail address:}
zxliu@email.jlu.edu.cn
(Zhenxin Liu).}\\
{\small College of Mathematics, Jilin University, Changchun 130012,
People's Republic of China}}
\date{}
\maketitle {\noindent{\bf Abstract}}

 The well-known Conley's theorem states that the complement of chain
recurrent set equals the union of all connecting orbits of the flow
$\varphi$ on the compact metric space $X$, i.e.
$X-\mathcal{CR}(\varphi)=\bigcup [B(A)-A]$, where
$\mathcal{CR}(\varphi)$ denotes the chain recurrent set of
$\varphi$, $A$ stands for an attractor and $B(A)$ is the basin
determined by $A$. In this paper we show that by appropriately
selecting the definition of random attractor, in fact we define a
random local attractor to be the $\omega$-limit set of some random
pre-attractor surrounding it, and by considering appropriate
measurability, in fact we also consider the universal
$\sigma$-algebra $\mathcal F^u$-measurability besides $\mathcal
F$-measurability, we are able to obtain the random case of Conley's
theorem.

\indent

\vskip 5mm

%\ams{37H99, 11B37, 37B35, 37B20, 37B25}

\section{Introduction and main result}

Among the tasks of differential equations and  dynamical systems, a
fundamental one is to study qualitative, asymptotic, long-term
behavior of solutions/orbits. Conley in his famous work \cite{Con}
introduced the concept of chain recurrence, and defined an attractor
to be the $\omega$-limit set of one of  its neighbourhoods. He
obtained the very interesting intrinsic relation between attractors
and chain recurrent set. First we take a simple retrospect about his
result.

Suppose $(X,d)$ is a compact metric space and $\varphi$ is a flow
with the phase space $X$. An open nonempty set $U$ is called a {\em
pre-attractor} for flow $\varphi$ if
\begin{equation}\label{pr}
\overline{\bigcup_{t\ge T}\varphi(t,U)}\subset U
\end{equation}
for some $T>0$. In fact, $\overline{\varphi(T_0,U)}\subset U$ for
some $T_0$ implies (\ref{pr}) holds, for details see page 33 of
\cite{Con}. For the convenience of late use, we adopt the form
(\ref{pr}). The {\em attractor} determined by the pre-attractor $U$
is defined by
\begin{equation}\label{DA}
A:=\bigcap_{t\ge T}\overline{\bigcup_{s\ge t}\varphi(s,U)}.
\end{equation}
It is easy to see that $A$ is a compact set, which is invariant
under the flow $\varphi$, i.e. $\varphi(t,A)=A,~\forall~ t\in
\mathbb R$. The {\em basin} of $A$, denoted by $B(A)$, is defined by
\[B(A)=\{x|~\varphi(t,x)\in U~{\rm for~some}~t\ge 0\}.\]
Since $X$ is compact, it is obvious that $B(A)$ is independent of
the choice of $U$. Therefore we denote it by $B(A)$, not mentioning
$U$.

For given $\epsilon, T>0$, a finite sequence $(x_0,t_0)$,
$(x_1,t_1)$, $\cdots$, $(x_n,t_n)$ in $X\times (0,\infty)$ is called
an $\epsilon$-$T$-chain for $\varphi$ if
\[ d(\varphi(t_j,x_j),x_{j+1})<\epsilon,~t_j\ge T
\]
for $j=0,1,\cdots,n-1$. And we call $n$ the {\em length} of the
chain. A point $p\in X$ is called {\em chain recurrent} if for any
$\epsilon, T>0$, there is an $\epsilon$-$T$-chain with the length at
least $1$ which begins and ends at $p$. And we use
$\mathcal{CR}(\varphi)$ to denote the  set of all chain recurrent
points in $X$.

Conley's theorem tells us that the complement of the chain recurrent
set is in fact the union of the sets $B(A)-A$, as $A$ varies over
the collection of attractors of $\varphi$, i.e.
\begin{equation}\label{con}
X-\mathcal{CR}(\varphi)=\bigcup [B(A)-A].
\end{equation}

For example, consider the differential equation $\dot{x}=x^3-x$ on
the interval $X=[-1,1]$ and assume $\varphi$ is the flow generated
by it. Besides the two trivial attractors $\emptyset,X$, it is
obvious that all other attractors are $\{0\}$, $[-1,0]$ and $[0,1]$
with basins of attractions $(-1,1)$, $[-1,1)$ and $(-1,1]$
respectively. Hence by (\ref{con}), we easily obtain that
$\mathcal{CR}(\varphi)=\{-1,0,1\}$, just the equilibria of the flow
$\varphi$. See the left picture in Figure \ref{fig}.

Another simple example, consider $\dot{\theta}=\cos^2\frac\theta2$
on $S^1$, and assume $\varphi$ is the flow generated by it. Then it
is easy to see that there is no nontrivial attractors for $\varphi$,
noticing that the unique equilibrium $\pi$ is not an attractor. So,
by (\ref{con}), the chain recurrent set for $\varphi$ is the whole
circle $S^1$. See the right picture in Figure \ref{fig}.
\begin{center}
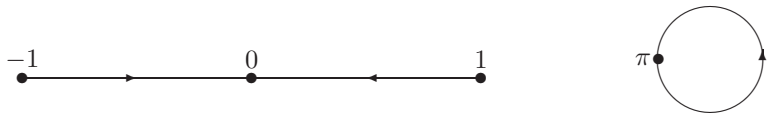
\begin{figure}\label{fig}
\setlength{\unitlength}{0.1in}
\begin{picture}(50,10)
\put(6,5){\circle*{0.5}} \put(6,5){\vector(1,0){6}}
\put(12,5){\line(1,0){6}} \put(18,5){\line(1,0){6}}
\put(18,5){\circle*{0.5}} \put(24,5){\line(1,0){6}}
\put(30,5){\vector(-1,0){6}} \put(30,5){\circle*{0.5}}
\put(42,6){\circle{8}} \put(39.3,6){\circle*{0.5}}
\put(44.75,6.0){\vector(0,1){0.5}} \put(6,6){\makebox(0,0){$-1$}}
\put(18,6){\makebox(0,0){$0$}}\put(30,6){\makebox(0,0){$1$}}
\put(38.5,6){\makebox(0,0){$\pi$}}
\end{picture}\caption{Chain recurrent sets for $\dot{x}=x^3-x$
and $\dot{\theta}=\cos^2\frac\theta2$}
\end{figure}
\end{center}

Conley's result was adapted for maps on compact spaces by Franks
\cite{Fra}, was later established for maps on locally compact metric
spaces by Hurley \cite{Hu0,Hu1}, and was extended by Hurley
\cite{Hu2} for semiflows and maps on arbitrary metric spaces. In
this paper, we will extend Conley's theorem to the random case, i.e.
we will show that the similar result holds for cocycle $\varphi$ on
compact metric spaces. For random dynamical systems, by defining
random chain recurrent variable, which is the counterpart of chain
recurrent point in random case, and by defining random attractor
similar to (\ref{DA}), we obtain the main result of this paper,
which states as follows:

\begin{theorem}{\rm\bf (Random Conley's theorem)}.
Suppose $U(\omega)$ is an arbitrary random pre-attractor,
$A(\omega)$ is the random local attractor determined by $U(\omega)$,
and $B(A)(\omega)$ is the random basin determined by $A(\omega)$,
then we have the following holds:
\begin{equation}\label{mr}
X-\mathcal{CR}_\varphi(\omega)=\bigcup [B(A)(\omega)-A(\omega)],
\end{equation}
where the union is taken over all random local attractors determined
by random pre-attractors, and ``=" in (\ref{mr}) holds $\mathbb
P$-almost surely.
\end{theorem}
Detailed definitions and notations in the main theorem can be found
in the next section.

Similar to deterministic Conley's theorem, our result accurately
describes where on earth the random chain recurrent variables lie.

\section{Preliminaries}
\setcounter{equation}{0}

In this section, we will give some preliminary definitions and
propositions for the late use. Firstly we give the definition of
continuous random dynamical systems (cf. Arnold \cite{Ar1}).
\begin{definition}
%{\rm (Random Dynamical System (RDS))}.
A (continuous) {\em random dynamical system (RDS)}, shortly denoted
by $\varphi$, consists
of two ingredients: \\
(i) A model of the noise, namely a metric dynamical system $(\Omega,
\mathcal F, \mathbb P, (\theta_t)_{t\in \mathbb T})$, where
$(\Omega, \mathcal F, \mathbb P)$ is a probability space and $(t,
w)\mapsto \theta_t\omega$ is a measurable flow which leaves $\mathbb
P$ invariant, i.e. $\theta_t\mathbb P=\mathbb P$ for all
$t\in \mathbb T$.\\
%For simplicity we also assume that $\theta$ is
%ergodic under $\mathbb P$, meaning that a $\theta$-invariant set
%has probability 0 or 1.
(ii) A model of the system perturbed by noise, namely a cocycle
$\varphi$ over $\theta$, i.e. a measurable mapping $\varphi: \mathbb
T\times \Omega\times X \rightarrow X,
(t,\omega,x)\mapsto\varphi(t,\omega,x)$, such that
$(t,x)\mapsto\varphi(t,\omega,x)$ is continuous for all
$\omega\in\Omega$ and the family
$\varphi(t,\omega,\cdot)=\varphi(t,\omega):X\rightarrow X$ of random
self-mappings of $X$ satisfies the cocycle property:
\begin{equation}\label{phi}
\varphi(0,\omega)={\rm id}_X, \varphi(t+s,\omega)=\varphi(t,\theta_s
\omega)\circ\varphi(s,\omega)\quad {\rm for ~all}\quad t,s\in\mathbb
T,\omega\in\Omega.
\end{equation}
In this definition, $\mathbb T=\mathbb Z$ or $\mathbb R$.
\end{definition}

It follows from (\ref{phi}) that $\varphi(t,\omega)$ is a
homeomorphism of $X$, and the fact
\[\varphi(t,\omega)^{-1}=\varphi(-t,\theta_t\omega)\]
is very useful in the following.

Below any mapping from $\Omega$ into the collection of all subsets
of $X$ is said to be a {\it multifunction} (or a set valued mapping)
from $\Omega$ into X. We now give the definition of random set,
which is a fundamental concept for RDS.

\begin{definition}
%{\rm\bf(Random Set)}.
Let $X$ be a metric space with a metric $d$. The multifunction
$\omega\mapsto D(\omega)\neq \emptyset$ is said to be a {\em random
set} if the mapping $\omega\mapsto {\rm dist}_X(x,D(\omega))$ is
measurable for any $x \in X$, where ${\rm dist}_X(x,B)$ is the
distance in $X$ between the element $x$ and the set $B \subset X$.
If $D(\omega)$ is closed/compact for each $\omega\in\Omega$,
$D(\omega)$ is called a {\em random closed/compact set}.
\end{definition}

Afterwards, we also call a multifunction $D(\omega)$ measurable for
convenience if the mapping $\omega\mapsto {\rm dist}_X(x,D(\omega))$
is measurable for any $x \in X$.

Now we enumerate some basic results about random sets in the
following propositions, for details the reader can refer to Arnold
\cite{Ar1}, Chueshov \cite{Chu} for instance.
\begin{proposition}
Let X be a Polish space, i.e. a separable complete metric space.
The following assertions hold: \\
(i) $D$ is a random set in $X$ if and only if the set $\{\omega :
D(\omega) \bigcap U \neq \emptyset \}$ is measurable for any open
set $U \subset X$;\\
(ii) D is a random set in $X$ if and only if $\overline{D(\omega)}$
is a random closed set ($\overline{D(\omega)}$ denotes the closure
of $D(\omega)$ in
$X$);\\
(iii) if D is a random open set, then the closure $\overline D$ of
$D$ is a random closed set; if D is a random closed set, then {\rm
int}$(D)$, the interior
of $D$, is a random open set; \\
(iv) $D$ is a random compact set in $X$ if and only if $D(\omega)$
is compact for every $\omega\in\Omega$ and the set $\{\omega :
D(\omega)\bigcap C \neq\emptyset\}$ is measurable for
any closed set $C \subset X$;\\
(v) if $\{D_n, n \in\mathbb N\}$ is a sequence of random closed sets
with non-void intersection and there exists $n_0 \in\mathbb N$ such
that $D_{n_0}$ is a random compact set, then
$\bigcap_{n\in\mathbb N}D_n$ is a random compact set in $X$;\\
(vi) if $\{D_n, n \in\mathbb N\}$ is a sequence of random sets, then
$D = \bigcup_{n\in\mathbb N}D_n$ is also a random set in $X$;
\end{proposition}
\begin{proposition} {\rm(Measurable Selection
Theorem)}. Let a multifunction $\omega\mapsto D(\omega)$ take values
in the subspace of closed non-void subsets of a Polish space $X$.
Then $D(\omega)$ is a random closed set if and only if there exists
a sequence $\{v_n : n \in\mathbb N\}$ of measurable maps $v_n
:\Omega\mapsto X$ such that \[ v_n(\omega)\in D(\omega)\quad and
\quad D(\omega) =\overline{\{v_n(\omega),n\in\mathbb N\}}\quad
for~all \quad \omega\in\Omega.\] In particular if $D(\omega)$ is a
random closed set, then there exists a measurable selection, i.e. a
measurable map $v :\Omega\mapsto X$  such that $v(\omega)\in
D(\omega)$ for all $\omega\in\Omega$.
\end{proposition}

Similar to deterministic case, we can define random chain
recurrence. The following random chain recurrent variable for random
dynamical systems is the counterpart of chain recurrent point for
the deterministic dynamical systems. Here, `recurrence' is defined
in the `pull-back' sense.
\begin{definition}\label{CR}
For given random variable $\epsilon(\omega)>0$, the $n+1$ pairs
$(x_0(\omega),t_0)$, $(x_1(\omega), t_1)$, $\cdots$, $(x_n(\omega),
t_n)$, where $x_0(\omega)$, $x_1(\omega)$, $\cdots$, $x_n(\omega)$
are random variables, are called a {\em random}
$\epsilon(\omega)$-$T(\omega)$-{\em chain}, if the following holds:
\[
d(\varphi(t_i,\theta_{-t_i}\omega)x_i(\theta_{-t_i}\omega),x_{i+1}(\omega))
<\epsilon(\omega),\quad i=0,\cdots,n-1,
\]
where $t_i\ge T(\omega)$, $T(\omega)$ is a positive random variable
almost surely. And we call $n$ the {\em length} of the random
$\epsilon(\omega)$-$T(\omega)$-chain. A random variable $x(\omega)$
is called {\em random chain recurrent} if for any given
$\epsilon(\omega),T(\omega)>0$, there exists an
$\epsilon(\omega)$-$T(\omega)$-chain beginning and ending at
$x(\omega)$ $\mathbb P$-almost surely; $x(\omega)$ is called {\em
partly random chain recurrent with index $\delta$} if for any
$\epsilon(\omega),T(\omega)>0$, there exists an
$\epsilon(\omega)$-$T(\omega)$-chain beginning and ending at
$x(\omega)$ with probability not less than $\delta$, where the index
$\delta$ is the maximal number satisfying this property; $x(\omega)$
is called {\em completely random non-chain recurrent} if there
exists $\epsilon_0(\omega),T_0(\omega)>0$ such that there is no
$\epsilon_0(\omega)$-$T_0(\omega)$-chain beginning and ending at
$x(\omega)$ with positive probability.
\end{definition}

\begin{remark}\rm
In the definition of partly random chain recurrence, the index
$\delta$ being maximal means that for
$\forall\eta>0,\exists\epsilon_0(\omega),T_0(\omega)>0$ such that
any $\epsilon_0(\omega)$-$T_0(\omega)$-chain begins and ends at
$x(\omega)$ with probability $\le\delta+\eta$.
\end{remark}

\begin{remark}\rm
In this paper, we will denote $\mathcal {CR}_\varphi(\omega)$ the
{\em random chain recurrent set of $\varphi$}, which has the
property that for any random chain recurrent variable $x(\omega)$,
we have $x(\omega)\in\mathcal{CR}_\varphi(\omega)$ $\mathbb
P$-almost surely, and vice versa (i.e. if a random variable
$x(\omega)\in\mathcal{CR}_\varphi(\omega)$ $\mathbb P$-almost
surely, then $x(\omega)$ is random chain recurrent); for any
completely random non-chain recurrent variable $x(\omega)$, we have
$x(\omega)\in X-\mathcal{CR}_\varphi(\omega)$ $\mathbb P$-almost
surely, and vice versa; for any partly random chain recurrent
variable $x(\omega)$ with index $\delta$, we have
$x(\omega)\in\mathcal{CR}_\varphi(\omega)$ with probability
$\delta$, and vice versa. For any given random variable $x(\omega)$,
denote
\[
\Omega_{\mathcal{CR}}(x)=\{\omega|~x(\omega)\in\mathcal
{CR}_\varphi(\omega)\}.
\]
If $x(\omega)$ is a partly random chain recurrent variable with
index $\delta$, then we call
$\{x(\omega)|~\omega\in\Omega_{\mathcal{CR}}(x)\}$ the {\em chain
recurrent  part of $x(\omega)$}. Therefore by the property of
$\mathcal{CR}_\varphi(\omega)$, we have that
$\mathcal{CR}_\varphi(\omega)$ is the union of all random chain
recurrent variables and the chain recurrent part of those partly
random chain recurrent variables.
\end{remark}

\begin{example}\rm
A simple example of random chain recurrent variable is random
equilibrium (a random variable $x(\omega)$ is called an equilibrium
if $\varphi(t,\omega)x(\omega)=x(\theta_t\omega)$ holds for $\forall
t>0,\omega\in\Omega$).
\end{example}

Throughout the paper, we will assume that $X$ is a compact metric
space, therefore it is a Polish space. The $\sigma$-algebra
$\mathcal F^u$ of universally measurable sets associated with the
base space $(\Omega,\mathcal F)$ is defined by $\mathcal
F^u=\bigcap_{\nu}{\bar{\mathcal F}}^\nu$, where the intersection is
taken over all probability measures $\nu$ on $(\Omega,\mathcal F)$
and ${\bar{\mathcal F}}^\nu$ denotes the completion of the
$\sigma$-algebra $\mathcal F$ with respect to the measure $\nu$. And
we call $\mathcal F^u$ the universal $\sigma$-algebra for brevity.
For a random variable $T(\omega)$, we call $T(\omega)>0$ if it holds
almost surely. By the measurable selection theorem, for any non-void
random closed set, there exists random variables belonging to it. In
the following, for a random open set, say $U(\omega)$, when we say
that a random variable $x(\omega)\in U(\omega)$, we mean that there
exists a random closed set $K(\omega)\subset U(\omega)$ such that
$x(\omega)\in K(\omega)$ almost surely.

For late use, we give the following important projection theorem,
which comes from \cite{Cas}.
\begin{proposition}{\rm (Projection Theorem)}.
Let X be a Polish space and $M\subset \Omega\times X$ be a set which
is measurable with respect to the product $\sigma$-algebra $\mathcal
F\times\mathcal B (X)$. Then the set
\[\Pi_\Omega M=\{\omega\in\Omega: (\omega,x)\in M ~{\rm for~ some}~ x\in X \}
\]
 is universally
measurable, i.e. belongs to $\mathcal F^u$, where $\Pi_\Omega$
stands for the canonical projection of ~$\Omega\times X$ to
$\Omega$. In particular it is measurable with respect to the
$\mathbb P$-completion ${\bar{\mathcal F}}^{\mathbb P}$  of
$\mathcal F$.
\end{proposition}

\begin{remark}\label{Fu}\rm
We have the following direct result. If $M\in{\mathcal
F}^u\times\mathcal B (X)$, $\Pi_\Omega M\in\mathcal F^u$ too. In
fact, we only need to show $(\mathcal F^u)^u=\mathcal F^u$ by
projection theorem. To see this, we notice that, on the one hand,
for arbitrary probability measure $\nu$ on measurable space
$(\Omega, \mathcal F^u)$, $\nu\mid_{\mathcal F}$, the restriction of
$\nu$ on $(\Omega, \mathcal F)$, is a probability measure on
$(\Omega, \mathcal F)$. On the other hand, for arbitrary probability
measure $\nu$ on $(\Omega, \mathcal F)$, we  can convert it into a
probability measure on $(\Omega, \mathcal F^u)$ by adding subsets of
$\Omega$ which are in $\mathcal F^u$ but not in $\mathcal F$ and
defining their measures to be $0$. That is, the measures on
$(\Omega, \mathcal F)$ and those on $(\Omega, \mathcal F^u)$ are one
to one. So by the fact that $\mathcal F\subset\mathcal
F^u\subset{\bar{\mathcal F}}^\nu$, where $\nu$ is an arbitrary
probability on $(\Omega, \mathcal F ~{\rm or}~\mathcal F^u)$, we
have
\[ {\bar{\mathcal F}}^\nu\subset(\overline{\mathcal F^u})^\nu\subset(\overline{\bar{\mathcal
F}^\nu})^\nu={\bar{\mathcal F}}^\nu.
\]
Therefore
\[\mathcal F^u=\bigcap_\nu{\bar{\mathcal F}}^\nu\subset
\bigcap_\nu(\overline{\mathcal F^u})^\nu=(\mathcal
F^u)^u\subset\bigcap_\nu{\bar{\mathcal F}}^\nu=\mathcal F^u,
\]
i.e. $(\mathcal F^u)^u=\mathcal F^u$.
\end{remark}

By remark \ref{Fu}, without loss of generality, we need only
consider $\mathcal F^u$-measurability throughout the paper.

At last, we give the definition of random local attractor and the
random basin determined by it.

\begin{definition}\label{UA}
%{\rm\bf(Random local attractor and Random basin)}.
A random open set $U(\omega)$ is called {\em random pre-attractor}
if it satisfies
\begin{equation}\label{at}
\overline{\bigcup_{t\ge
T(\omega)}\varphi(t,\theta_{-t}\omega)U(\theta_{-t}\omega)} \subset
U(\omega)\quad for~ some \quad T(\omega)>0,
\end{equation}
where $T(\omega)$ is an $\mathcal F$-measurable random variable. And
we define the {\em random local attractor} $A(\omega)$ inside
$U(\omega)$ to be the following:
\[A(\omega)=\bigcap_{n\in\mathbb N}\overline {\bigcup_{s\ge
nT(\omega)}\varphi(s,\theta_{-s}\omega)U(\theta_{-s}\omega)}.
\]
And the {\em random basin} $B(A)(\omega)$ determined by $A(\omega)$
is defined as follows
\begin{equation}\label{BA}
B(A)(\omega)=\{x: \varphi(t,\omega)x\in U(\theta_t\omega) \quad {\rm
for~ some}~ t\ge 0\}.
\end{equation}
\end{definition}

It is easy to see that in the above definition,  the random basin
$B(A)(\omega)$ may depend on the pre-attractor $U(\omega)$. In fact,
we can show that the basin is independent of the choice of
$U(\omega)$ and we defer the proof to the next section.

\begin{example}\rm
Suppose $x(\omega)$ is a random variable, and
$\epsilon_0(\omega),T_0(\omega)$ are two positive  random variables.
Consider the random set $U_x(\omega)$ determined by (\ref{U}) in the
next section. By the proof of lemma \ref{3.7}, we know that
$U_x(\omega)$ is a random pre-attractor, but it is not necessarily a
forward invariant random set. In fact, for $t<T_0(\omega)$, we can
not obtain that (\ref{pre}) also holds, which guarantees that
$U_x(\omega)$ is forward invariant. Hence the attractor
$A_x(\omega)$ determined by $U_x(\omega)$ is not necessarily the
attractor defined in \cite{Cra4}.
\end{example}

\section{Proof of the main result}
\setcounter{equation}{0}

Denote $U(T(\omega))=\overline {\bigcup_{s\ge
T(\omega)}\varphi(s,\theta_{-s}\omega)U(\theta_{-s}\omega)}$.
\begin{lemma}\label{A}
Suppose $U(\omega)$ is a given pre-attractor, then $U(T(\omega))$
and the the random local attractor
\begin{equation}\label{a}
A(\omega)=\bigcap_{n\in\mathbb N}U(nT(\omega))
\end{equation}
determined by $U(\omega)$ are random closed sets measurable with
respect to $\mathcal F^u$. Moreover, $A(\omega)$ is invariant, i.e.
$\varphi(t,\omega)A(\omega)=A(\theta_t\omega)$ for all $t\ge 0$, and
$A(\omega)$ is a local random pull-back set attractor, therefore a
local weak random set attractor in $U(\omega)$ (for the definition
of weak random set attractor see \cite{Och}).
\end{lemma}

\noindent{\bf Proof.} (i) We first show that $U(T(\omega))$ is a
random closed set. The idea of the proof is borrowed from
\cite{Cra5} and \cite{Chu}. For every $x\in X$, define
\begin{equation}\label{dtw}
d(t,\omega):={\rm
dist}_X(x,\varphi(t,\theta_{-t}\omega)U(\theta_{-t}\omega)).
\end{equation}
By (ii) of proposition 2.1 and the proof of proposition 1.5.1 of
\cite{Chu}, we obtain that the function $(t,\omega)\mapsto
d(t,\omega)$ is $\mathcal B(\mathbb T)\times\mathcal F$-measurable.
Clearly we have
\begin{eqnarray}
{\rm dist}_X(x,U(T(\omega)))&={\rm dist}_X(x,\bigcup_{t\ge
T(\omega)}\varphi(t,\theta_{-t}\omega)U(\theta_{-t}\omega))\nonumber\\
 &=
{\rm inf}_{t\ge T(\omega)}d(t,\omega).
\end{eqnarray}

For arbitrary $a\in\mathbb R^+$, we have
\[
\{\omega|~{\rm inf}_{t\ge
T(\omega)}d(t,\omega)<a\}=\Pi_\Omega\{(t,\omega)|~d(t,\omega)<a,t\ge
T(\omega)\}.
\]
It is obvious that the function $(t,\omega)\mapsto t-T(\omega)$ is
measurable with respect to $\mathcal B(\mathbb T)\times\mathcal F$,
so by projection theorem, we obtain that $\{\omega|~{\rm inf}_{t\ge
T(\omega)}d(t,\omega)<a\}$ is $\mathcal F^u$-measurable, which
follows that $U(T(\omega))$ is a random set measurable with respect
to the universal $\sigma$-algebra $\mathcal F^u$. The closeness of
$U(T(\omega))$ is obvious.

(ii) Clearly we have $U(nT(\omega))\supset U((n+1)T(\omega))$, which
follows that
\begin{equation}\label{aa}
A(\omega)={\rm lim}_{n\rightarrow\infty}U(nT(\omega)).
\end{equation}
Therefore $A(\omega)$ is  a random closed set. To get the attraction
property of $A(\omega)$, we notice that, for any given random
compact set $K(\omega)\subset U(\omega)$, \[{\rm
lim}_{t\rightarrow\infty}
d(\varphi(t,\theta_{-t}\omega)K(\theta_{-t}\omega)|A(\omega))=0
\] holds almost surely by (\ref{aa}), where the metric $d(A|B)$ between
two closed sets stands for the Hausdorff semi-metric, i.e.
$d(A|B):={\rm sup}_{x\in A}{\rm inf}_{y\in B}d(x,y)$. Hence
$A(\omega)$ is a pull-back set attractor in $U(\omega)$. Then by the
property of measure preserving of $\{\theta_t\}_{t\in\mathbb T}$, we
obtain that
\[
\mathbb P-{\rm lim}_{t\rightarrow\infty}{\rm d}
(\varphi(t,\omega)K(\omega),A(\theta_t\omega))=0,
\]
which implies that $A(\omega)$ is a weak random set attractor in
$U(\omega)$.

The rest work is to verify the invariance of $A(\omega)$. The
forward invariance of $A(\omega)$ follows from \cite{Cra1}, just
changing a few details. For completeness, we give its proof here.
For arbitrary $t\ge 0$,
\begin{align*}
\varphi(t,\omega)A(\omega)&= \varphi(t,\omega)\bigcap_{n\in\mathbb
N}\overline{\bigcup_{\tau\ge
nT(\omega)}\varphi(\tau,\theta_{-\tau}\omega)U(\theta_{-\tau}\omega)}\nonumber\\
&\subset\bigcap_{n\in\mathbb
N}\varphi(t,\omega)\overline{\bigcup_{\tau\ge
nT(\omega)}\varphi(\tau,\theta_{-\tau}\omega)U(\theta_{-\tau}\omega)}\nonumber\\
&\subset\bigcap_{n\in\mathbb N}\overline{\bigcup_{\tau\ge
nT(\omega)}\varphi(t,\omega)\varphi(\tau,\theta_{-\tau}\omega)U(\theta_{-\tau}\omega)}
\nonumber\\
&=\bigcap_{n\in\mathbb N}\overline{\bigcup_{\tau\ge
nT(\omega)}\varphi(t+\tau,\theta_{-(t+\tau)}\circ\theta_t\omega)
U(\theta_{-(t+\tau)}\circ\theta_t\omega)}\nonumber\\
&=\bigcap_{n\in\mathbb N}\overline{\bigcup_{\tau\ge
nT(\omega)+t}\varphi(\tau,\theta_{-\tau}\circ\theta_t\omega)
U(\theta_{-\tau}\circ\theta_t\omega)}\nonumber\\
&\subset\bigcap_{n\in\mathbb N}\overline{\bigcup_{\tau\ge
nT(\omega)}\varphi(\tau,\theta_{-\tau}\circ\theta_t\omega)
U(\theta_{-\tau}\circ\theta_t\omega)}=A(\theta_t\omega),\nonumber
\end{align*}
where the first two inclusions follows from the facts
$f(\bigcap_\alpha A_\alpha)\subset\bigcap_\alpha f(A_\alpha)$ for
arbitrary $f$ and $f(\bar A)\subset\overline{f(A)}$ for $f$
continuous respectively. The backward invariance of $A(\omega)$ is
similar to \cite{Cra5} noting that $X$ is compact, so we omit the
details here. This completes the proof of the lemma. \hfill$\square$

\begin{lemma}
Suppose $A(\omega)$ is a random local attractor and
$U_1(\omega),U_2(\omega)$ are two pre-attractors which determine the
same attractor $A(\omega)$, then the two basins determined by
$U_1,U_2$ respectively are equal $\mathbb P$-almost surely.
\end{lemma}

\noindent{\bf Proof.} Denote $B_1(A)(\omega),B_2(A)(\omega)$ the
basins determined by $U_1(\omega),U_2(\omega)$ respectively.  For
arbitrary random variable $x(\omega)\in B_1(A)(\omega)$, there
exists $t(\omega)\ge 0$ such that
\[\varphi(t(\omega),\omega)x(\omega)\in U_1(\theta_{t(\omega)}\omega),~
\forall ~\omega\in \Omega\] by the definition of basin. By
attraction property of $A(\omega)$ and the measure preserving  of
$\theta_t$, it follows that for $\forall~ \epsilon >0$, the
following holds:
\begin{align}\label{B}
&{\rm lim}_{s\rightarrow \infty}\mathbb P\{\omega|~{\rm
dist}_X(\varphi(s,\theta_{t(\omega)}\omega)\varphi(t(\omega),\omega)x(\omega),
A(\theta_s\circ\theta_{t(\omega)}\omega))>\epsilon\}\nonumber
\\
=&{\rm lim}_{s\rightarrow \infty}\mathbb P\{\omega|~{\rm
dist}_X(\varphi(s+t(\omega),\omega)x(\omega),
A(\theta_{s+t(\omega)}\omega))>\epsilon\}\nonumber
\\
=&{\rm lim}_{s\rightarrow \infty}\mathbb P\{\omega|~{\rm
dist}_X(\varphi(s,\omega)x(\omega),
A(\theta_s\omega))>\epsilon\}\nonumber
\\
=&{\rm lim}_{s\rightarrow \infty}\mathbb P\{\omega|~{\rm
dist}_X(\varphi(s,\theta_{-s}\omega)x(\theta_{-s}\omega),
A(\omega))>\epsilon\}=0.
\end{align}
Denote $d(\omega):={\rm dist}_X(A(\omega),X-U_2(\omega))$, where the
distance ${\rm dist}_X(A,B)$ between two closed sets stands for the
smallest distance between them, i.e. ${\rm dist}_X(A,B)={\rm
inf}_{x\in A}{\rm inf}_{y\in B}d(x,y)$. Therefore we have
$d(\omega)>0$ almost surely by the compactness of $X$. By a standard
argument, we obtain that there exists a $\delta>0$ such that
\begin{equation}
\mathbb P\{\omega|~ d(\omega)>\delta\}\ge 1-\epsilon.
\end{equation}
Denote $\Omega_\delta:=\{\omega|~d(\omega)>\delta\}$. By (\ref{B})
we have
\[
{\rm lim}_{s\rightarrow\infty}\mathbb P\{\omega\in\Omega_\delta|~
d(\varphi(s,\theta_{-s}\omega
)x(\theta_{-s}\omega),A(\omega))>\frac\delta 4\}=0.
\]
Therefore there exists $t_\delta(\omega)\ge 0$ such that
\[
d(\varphi(s,\theta_{-s}\omega)x(\theta_{-s}\omega),A(\omega))<\frac\delta
2,~ s\ge t_\delta(\omega)
\]
holds almost surely on $\Omega_\delta$. Hence by the definition of
$B_2(A)(\omega)$ and the measure preserving of $\theta_t$, we get
that
\begin{eqnarray}
\mathbb P\{\omega|~x(\omega)\in B_2(A)(\omega)\}&=\mathbb
P\{\omega|~\varphi(s,\theta_{-s}\omega)x(\theta_{-s}\omega)\in
U_2(\omega)~{\rm for~some}~s\ge 0\}\nonumber\\
&\ge \mathbb P(\Omega_\delta)\ge 1-\epsilon.\nonumber
\end{eqnarray}
Since $\epsilon>0$ is arbitrary, we have $x(\omega)\in
B_2(A)(\omega)$ almost surely by letting $\epsilon\rightarrow 0$. It
follows that $B_1(A)(\omega)\subset B_2(A)(\omega)$ almost surely,
and the converse inclusion is similar. This terminates the proof of
the lemma.\hfill$\square$

The above lemma indicates that the basin $B(A)(\omega)$ is uniquely
determined by $A(\omega)$, therefore is well defined. The following
lemma says that the basin is backward invariant random open.

\begin{lemma}\label{Ba}
For any given random local attractor $A(\omega)$, the random basin
$B(A)(\omega)$ determined by $A(\omega)$ is a backward invariant
random open set.
\end{lemma}

\noindent{\bf Proof.} It is obvious that $x\in B(A)(\omega)$ is
equivalent to $\varphi(t,\omega)x\in U(\theta_t\omega)$ for some
 $t\ge 0$, i.e. $x\in
\varphi(-t,\theta_{t}\omega)U(\theta_{t}\omega)$ by using the fact
that $\varphi(t,\omega)^{-1}=\varphi(-t,\theta_{t}\omega)$. So it
follows that $x\in B(A)(\omega)$ if and only if $x\in \bigcup_{t\ge
0}\varphi(-t,\theta_{t}\omega)U(\theta_{t}\omega)$, therefore we
obtain that
\[
B(A)(\omega)=\bigcup_{t\ge
0}\varphi(-t,\theta_{t}\omega)U(\theta_{t}\omega).
\]
 Then by a
similar method to that of \cite{Cra5,Cra1,Chu}, we can easily obtain
that
\[\overline{\bigcup_{t\ge
0}\varphi(-t,\theta_{t}\omega)\overline{U(\theta_{t}\omega)}}\] is
$\mathcal F^u$-measurable, hence it is a random closed set. It is
easy to verify that we have the following holds
\[
\overline{\bigcup_{t\ge
0}\varphi(-t,\theta_{t}\omega)U(\theta_{t}\omega)}=
\overline{\bigcup_{t\ge
0}\varphi(-t,\theta_{t}\omega)\overline{U(\theta_{t}\omega)}}.
\]

To see this, note first that the left hand is obvious the subset of
the right hand. And every element of the union of right hand is a
subset of the left hand, so the closure of the union of right is
included by the left for the closeness of the left hand. So we have
got that the closure of $B(A)(\omega)$ is a random closed set. By
(ii) of proposition 2.1 we obtain that $B(A)(\omega)$ is a random
set, the openness of $B(A)(\omega)$ follows the fact that
$\varphi(-t,\theta_t\omega)$ is homeomorphism on $X$.

To show the backward invariance of $B(A)(\omega)$, we only need to
show that its complement $X-B(A)(\omega)$ is forward invariant. If
the assertion is false, then there exists an $x_0\in X-B(A)(\omega)$
and $t_0>0$ such that $\varphi(t_0,\omega)x_0\in
B(A)(\theta_{t_0}\omega)$. Then by the definition of random basin,
we have $\varphi(t_1,\theta_{t_0}\omega)\varphi(t_0,\omega)x_0\in
U(\theta_{t_0+t_1}\omega)$ for some $t_1\ge 0$, where $U$ is a
random pre-attractor which determines $A$. But by the definition of
random basin, $X-B(A)(\omega)$ is the set of points that never enter
$U$, a contradiction. This terminates the proof of the
lemma.\hfill$\square$

By the above preliminaries, we can prove our main theorem now.  We
decompose the proof of the main theorem into the following several
lemmas:

\begin{lemma}\label{lemu}
If the random chain recurrent variable $x(\omega)\in U(\omega)$
$\mathbb P$-almost surely, where $U(\omega)$ is a random
pre-attractor, then we have $x(\omega)\in A(\omega)$ $\mathbb
P$-almost surely, where $A(\omega)$ is the attractor determined by
$U(\omega)$.
\end{lemma}

\noindent{\bf Proof.} If $x(\omega)\in {U(\omega)}$ $\mathbb
P$-almost surely, recalling that $U(\omega )$ is a random
pre-attractor, fix $T(\omega)>0$ such that $U(T(\omega))\subset
U(\omega)$. Denote
  \begin{equation}\label{d}
    \bar{d}(\omega):={\rm dist}_X(U(T(\omega)), X-U(\omega)).
  \end{equation}

By the compactness of $X$, it is obvious that $\bar{d}(\omega)>0$
holds almost surely. By measurable selection theorem, similar to
\cite{Ar2}, there exists two collections of random variables
$\{x_n(\omega)\}_{n\in\mathbb N}\subset U(T(\omega))$ with
$U(T(\omega))=\overline{\{x_n(\omega)\}}_{n\in\mathbb N}$ and
$\{x_m(\omega)\}_{m\in\mathbb N}\subset X-U(\omega)$ with
$X-U(\omega)=\overline{\{x_m(\omega)\}}_{m\in\mathbb N}$ such that
\begin{equation}
\bar{d}(\omega)={\rm inf}_{n\in \mathbb N,m\in \mathbb
N}d(x_n(\omega),x_m(\omega)).
\end{equation}

So we obtain that $\bar{d}(\omega)$ is a random variable. Choose
$0<\epsilon(\omega)<\bar{d}(\omega)$, then we have
\[\varphi(t,\theta_{-t}\omega)x(\theta_{-t}\omega)\in
\varphi(t,\theta_{-t}\omega)U(\theta_{-t}\omega)\subset
U(T(\omega)), ~~ {\rm where} ~t\ge T(\omega).\]

Then by the the fact that $x(\omega)$ is random chain recurrent, for
this $\epsilon(\omega)
>0$ and the above $T(\omega)$, there exists an $\epsilon(\omega)$-$T(\omega)$-chain
$(x_0(\omega),t_0)$, $(x_1(\omega), t_1)$, $\cdots$, $(x_n(\omega),
t_n)$ with $x_0(\omega)=x_n(\omega)=x(\omega)$. So by the choice of
$\epsilon(\omega)$ and induction we obtain that
\begin{equation}\label{p}
x(\omega)\in B_{\epsilon(\omega)}(U(T(\omega)))\subset U(\omega)
\end{equation}
almost surely, where $B_r(x)$ stands for the open ball centered at
$x$ with radius $r$.

Hence we derive $x(\omega)\in U(T(\omega))$ from letting
$\epsilon(\omega)\rightarrow 0$ in (\ref{p}) and from the closeness
of $U(T(\omega))$. And then let $T(\omega)\rightarrow \infty$ in
(\ref{p}), we obtain $x(\omega)\in A(\omega)$ almost surely by
(\ref{aa}). \hfill$\square$

\begin{lemma}\label{tome}
Suppose $U(\omega)$ is a random open set, $x(\omega)$ is a random
variable. Define
\begin{equation}\label{t}
t(\omega):={\rm inf}\{t\in\mathbb R^+ |~
\varphi(t,\omega)x(\omega)\in \overline {U(\theta_t \omega)}\},
\end{equation}
i.e. the first entrance time of $x(\omega)$ into $U(\omega)$ under
the cocycle $\varphi$. Then $t(\omega)$ is a random variable, which
is measurable with respect to the universal $\sigma$-algebra
$\mathcal F^u$.
\end{lemma}
\noindent{\bf Proof.} It is easy to see that
\[
t(\omega)={\rm inf}\{t\in\mathbb R^+ |~ {\rm
dist}_X(\varphi(t,\omega)x(\omega),\overline {U(\theta_t \omega)})=0
\}.
\]
Since the function
\[
(t,\omega)\mapsto {\rm dist}_X(\varphi(t,\omega)x(\omega),\overline
{U(\theta_t \omega)})=:\tilde d(t,\omega)
\]
is $\mathcal B(\mathbb T)\times\mathcal F$-measurable by a similar
argument as the proof of lemma \ref{A}. For arbitrary $a\in \mathbb
R^+$, it is easy to see that
\[
\{\omega|~t(\omega)<a\}=\Pi_\Omega\{(t,\omega)|~0\le t<a, \tilde
d(t,\omega)=0\}.
\]
It is obvious that $\{(t,\omega)|~0\le t<a, \tilde
d(t,\omega)=0\}\in \mathcal B(\mathbb T)\times\mathcal F$, so
$\{\omega|~t(\omega)<a\}$ is $\mathcal F^u$-measurable by projection
theorem. \hfill$\square$

\begin{lemma}\label{lemb}
If the random chain recurrent variable $x(\omega)\in B(A)(\omega)$
$\mathbb P$-almost surely, then we have $x(\omega)\in A(\omega)$
$\mathbb P$-almost surely, where $B(A)(\omega)$ is the basin
determined by the random local attractor $A(\omega)$.
\end{lemma}
\noindent{\bf Proof.} For $x(\omega)\in B(A)(\omega)$, take
$t(\omega)\ge 0$ defined by (\ref{t}). Suppose $U(\omega)$ is a
pre-attractor which determines the attractor $A(\omega)$. So we
have, for $s\ge T(\theta_{t(\omega)}\omega)$, the following almost
surely holds:
\begin{align}
\varphi(s+t(\omega),\omega)x(\omega)&=\varphi(s,\theta_{t(\omega)}
\omega)\varphi(t(\omega),\omega)x(\omega)\nonumber\\
&\in\overline{\varphi(s,\theta_{t(\omega)}\omega)U(\theta_{t(\omega)}\omega)}\subset^\ast
U(\theta_{s+t(\omega)}\omega).\nonumber
\end{align}
In fact, by the definition of pre-attractor, there exists a random
variable $T(\theta_{t(\omega)}\omega)>0$ such that
\[
\overline{\bigcup_{s\ge
T(\theta_{t(\omega)}\omega)}\varphi(s,\theta_{-s}\circ\theta_{t(\omega)}\omega)
U(\theta_{-s}\circ\theta_{t(\omega)}\omega)}\subset
U(\theta_{t(\omega)}\omega).
\]
Then by the property of measure preserving of $\theta_t$, we obtain
that $\subset^*$ holds $\mathbb P$-almost surely. Hence we obtain
\begin{align}
&\mathbb P\{\omega|~\varphi(s+t(\omega),\omega)x(\omega)\in
U(\theta_{s+t(\omega)}\omega), s\ge T(\theta_{t(\omega)}\omega)\}\nonumber\\
=&\mathbb
P\{\omega|~\varphi(s+t(\omega),\theta_{-(s+t(\omega))}\omega)
x(\theta_{-(s+t(\omega))}\omega)\in U(\omega), s\ge T(\theta_{t(\omega)}\omega)\}\nonumber\\
=&1 \nonumber
\end{align}
by the fact that $\mathbb P$ is invariant under $\theta_t$ again.
Further more, we are able to obtain the following finer result:
\begin{align}
&\quad\varphi(s+t(\omega),\theta_{-(s+t(\omega))}\omega)
x(\theta_{-(s+t(\omega))}\omega)\nonumber\\
&=
\varphi(s,\theta_{-s}\omega)\varphi(t(\omega),\theta_{-(s+t(\omega))}\omega)
x(\theta_{-(s+t(\omega))}\omega)\nonumber\\
&\in^* \overline{\varphi(s,\theta_{-s}\omega)
U(\theta_{-s}\omega)}\subset^{**} U(\omega), \nonumber
\end{align}
where $\in^*$ holds using the property of measure preserving of
$\theta_t$ again. And $\subset^{**}$ holds if $s\ge T(\omega)$ by
the property of pre-attractor. So denote $\tilde T(\omega)={\rm
max}\{T(\omega),T(\theta_{t(\omega)}\omega)\}$, and take
$T_1(\omega)=t(\omega)+\tilde T(\omega)$, then it follows that
\[\overline{\varphi(s,\theta_{-s}\omega)
U(\theta_{-s}\omega)}\subset U(T(\omega))\subset U(\omega)\]
whenever $s\ge T_1(\omega)$. Take $\bar{d}(\omega)$ as defined by
(\ref{d}), fix $T_1(\omega)$, choose $0<\epsilon(\omega)<d(\omega)
$. Then it follows that any random
$\epsilon(\omega)$-$T_1(\omega)$-chain of length one begins at
$x(\omega)\in B(A)(\omega)$ must ends in
$B_{\epsilon(\omega)}(U(T(\omega)))$ almost surely. By the fact that
$x(\omega)$ is a random chain recurrent variable and repeat the
proof process of lemma \ref{lemu} we obtain that $x(\omega)\in
A(\omega)$ almost surely. \hfill$\square$

\begin{remark}\label{rem}\rm
By the process of proofs of lemma \ref{lemu} and lemma \ref{lemb},
if
\[
\mathbb P\{\omega|~x(\omega)\in B(A)(\omega)\}=\delta<1,~{\rm and}~
x(\omega)~{\rm is~random~chain~recurrent},
\]
then we can easily obtain that
\[
\mathbb P\{\omega|~x(\omega)\in A(\omega)\}=\delta,
\]
i.e.
\[
\mathbb P\{\omega|~x(\omega)\in B(A)(\omega)\backslash
A(\omega)\}=0.
\]
If $x(\omega)$ is partly random chain recurrent with index $\delta$,
then we have if $\omega\in\Omega_{\mathcal{CR}}(x)$, then
\[
x(\omega)\in A(\omega)~{\rm or}~x(\omega)\in X\backslash
B(A)(\omega)
\]
holds except for a null set for any given attractor $A(\omega)$ by
the proof process of lemma \ref{lemu} and lemma \ref{lemb} again.
Hence we obtain that
\[
\omega\in\Omega\backslash\Omega_{\mathcal{CR}}(x)~{\rm whenever}~
x(\omega)\in B(A)(\omega)-A(\omega),
\]
i.e.
\[
x(\omega)\in X-\mathcal{CR}_\varphi(\omega)~{\rm whenever}~
x(\omega)\in B(A)(\omega)-A(\omega).
\]
Therefore by lemma \ref{lemu}, lemma \ref{lemb} and above arguments
we obtain that
\[
B(A)(\omega)-A(\omega)\subset X-\mathcal{CR}_\varphi(\omega)
\]
holds $\mathbb P$-almost surely for any given attractor $A(\omega)$.
When the number of attractors is uncountable, to avoid the
possibility that
\[
\bigcup[B(A)(\omega)-A(\omega)]\subset
X-\mathcal{CR}_\varphi(\omega)
\]
does not hold $\mathbb P$-almost surely, we can redefine
$B(A)(\omega)-A(\omega)$ on the null set such that
\[
B(A)(\omega)-A(\omega)\subset X-\mathcal{CR}_\varphi(\omega)~{\rm\bf
for~all}~\omega\in\Omega.
\]
In fact, if $B(A)(\omega)-A(\omega)\subset
X-\mathcal{CR}_\varphi(\omega)$ holds on a set $\tilde{\Omega}$ of
full measure, we can set
\[
B(A)(\omega)-A(\omega)=X-\mathcal{CR}_\varphi(\omega)~{\rm
on}~\Omega\backslash\tilde{\Omega}.
\]
Hence by this way, we have obtained
\[
\bigcup[B(A)(\omega)-A(\omega)]\subset
X-\mathcal{CR}_\varphi(\omega)
\]
holds $\mathbb P$-almost surely.
\end{remark}

By lemma \ref{lemu}, lemma \ref{lemb}, and  remark \ref{rem} we
obtain that the right hand of (\ref{mr}) is a subset of the left
hand $\mathbb P$-almost surely. To prove the equality (\ref{mr}),
the rest work is to verify that the converse inclusion is true
$\mathbb P$-almost surely. To this end, we first define a random
open set $U_x(\omega)$ measurable with respect to $\mathcal F^u$ for
late use, which is defined as follows.

Suppose $x(\omega)$ is a random variable,
$\epsilon_0(\omega),T_0(\omega)$ are two positive  random variables.
Define
\begin{eqnarray}
U_1(\omega):&=\bigcup_{t\ge
T_0(\omega)}B_{\epsilon_0(\omega)}(\varphi(t,\theta_{-t}\omega)x(\theta_{-t}\omega)),
\nonumber\\
U_2(\omega):&=\bigcup_{t\ge
T_0(\omega)}B_{\epsilon_0(\omega)}(\varphi(t,\theta_{-t}\omega)U_1(\theta_{-t}\omega)),
\nonumber\\
 &\cdots\cdots\cdots \cdots\cdots                      \nonumber\\
U_n(\omega):&=\bigcup_{t\ge
T_0(\omega)}B_{\epsilon_0(\omega)}(\varphi(t,\theta_{-t}\omega)U_{n-1}(\theta_{-t}\omega)),
\nonumber\\
&\cdots\cdots\cdots \cdots\cdots~~. \nonumber
\end{eqnarray}
By the proof method of proposition 1.5.1 of \cite{Chu} on page 32,
which in turn stems from \cite{Cra5}, and a similar argument as in
lemma \ref{Ba}, adding some slight changes in the process of proof,
we can conclude that $U_1(\omega),\cdots,U_n(\omega),\cdots$ are all
$\mathcal F^u$-measurable open sets. We omit the details here. So
the set
\begin{equation}\label{U}
U_x(\omega):=\bigcup_{n\in \mathbb N}U_n(\omega)
\end{equation}
is a random open set measurable with respect to $\mathcal F^u$ by
(vi) of proposition 2.1.

Now we can give the proof of the converse inclusion of (\ref{mr}),
see the following lemma.

\begin{lemma}\label{3.7}
$X-\mathcal {CR}_\varphi(\omega)\subset\bigcup
[B(A)(\omega)-A(\omega)]$ holds $\mathbb P$-almost surely.
\end{lemma}
\noindent{\bf Proof.} We divide the proof into two cases.\\
Case 1: when $x(\omega)$ is completely random non-chain recurrent,
i.e. $x(\omega)\in X-\mathcal{CR}_\varphi(\omega)$ $\mathbb
P$-almost surely, there exists $\epsilon_0(\omega)>0, T_0(\omega)>0$
such that there is no $\epsilon_0(\omega)$-$T_0(\omega)$-chain which
begins and ends at $x(\omega)$ with positive probability. Take
$U_x(\omega)$ defined by (\ref{U}), then it is easy to see that
$x(\omega)\notin U_x(\omega)$ and
\begin{equation}\label{pre}
\varphi(t,\theta_{-t}\omega)x(\theta_{-t}\omega)\in U_x(\omega)
\end{equation}
 when $t\ge T_0(\omega)$. Given $y(\omega)\in U_x(\omega)$,
it is obvious that the $\epsilon_0(\omega)$-neighbourhood of
$\varphi(t,\theta_{-t}\omega)y(\theta_{-t}\omega)$ lies in
$U_x(\omega)$ whenever $t\ge T_0(\omega)$, therefore
${U_x(T_0(\omega))}\subset U_x(\omega)$, where $U_x(T_0(\omega))$ is
defined similarly to $U(T(\omega))$. So $U_x(\omega)$ is a random
pre-attractor and it determines a random local attractor
$A_x(\omega)$ by lemma \ref{A} (The only difference is that
$U_x(\omega)$ is $\mathcal F^u$-measurable while the $U(\omega)$ in
lemma \ref{A} is $\mathcal F$-measurable. But by remark \ref{Fu} we
know that this does not affect the result). And we have
$x(\omega)\in B(A_x)(\omega)-U_x(\omega)\subset
B(A_x)(\omega)-A_x(\omega)$ by (\ref{pre}). \\
Case 2: when $x(\omega)$ is partly random chain recurrent with index
$\delta$, i.e. $x(\omega)\in X-\mathcal{CR}_\varphi(\omega)$ with
probability $1-\delta$, then for $\forall\eta
>0,\exists\epsilon_0(\omega),T_0(\omega)>0$ such that any
$\epsilon_0(\omega)$-$T_0(\omega)$-chain begins and ends at
$x(\omega)$ with probability $\le\delta+\eta$. Take $U_x(\omega)$
defined by (\ref{U}), then by the proof of case 1, it is easy to see
that $x(\omega)\in B(A_x)(\omega)$ $\mathbb P$-almost surely and
\[
\mathbb P\{\omega|~x(\omega)\in U_x(\omega)\}\le \delta+\eta,
\]
hence
\[
\mathbb P\{\omega|~x(\omega)\in  B(A_x)(\omega)-A_x(\omega)\}\ge
1-(\delta+\eta).
\]
Thus we obtain
\[
\mathbb P\{\omega|~x(\omega)\in\bigcup [B(A)(\omega)-A(\omega)]\}\ge
1-\delta
\]
by letting $\eta\rightarrow 0$.

Therefore, by the proof of case 1 and case 2, we obtain
\[
X-\mathcal{CR}_\varphi(\omega)\subset\bigcup
[B(A)(\omega)-A(\omega)]
\]
holds $\mathbb P$-almost surely. This completes the proof of the
lemma. \hfill$\square$

By lemmas \ref{lemu}, \ref{lemb} and \ref{3.7}, and remark
\ref{rem}, we complete the proof of our main theorem.

\section{Two simple examples}

In this section, we will give two examples to illustrate our
results. The first example is borrowed from \cite{Cra4}.
\begin{example}\rm
Consider the Stratonovich stochastic differential equation (SDE)
\[
{\rm d}X_t=(X_t-X_t^3){\rm d}t+(X_t-X_t^3)\circ{\rm d}W_t
\]
on the interval $[-1,1]$. From p.123 of \cite{Klo} we know that the
RDS $\varphi:\mathbb R\times\Omega\times [-1,1]\mapsto [-1,1]$
generated by this SDE can be expressed by
\[
\varphi(t,\omega)x=\frac{x{\rm e}^{t+W_t(\omega)}}{(1-x^2+x^2{\rm
e}^{2t+2W_t(\omega)})^{\frac12}}.
\]
It is easily to see that $A_1=\{-1\}$, $A_2=\{1\}$, $A_3=\{-1,1\}$
are all attractors except for the two trivial attractors
$\emptyset,[-1,1]$, and the corresponding basins of attraction are
$[-1,0)$, $(0,1]$ and $[-1,0)\bigcup(0,1]$ respectively. Hence by
our main theorem, we obtain that the random chain recurrent
variables are $x_1(\omega)\equiv -1$, $x_2(\omega)\equiv 1$ and
$x_3(\omega)\equiv 0$ and the combinations of $-1,0,1$, e.g.
$x(\omega)=-1,0,1$ with probabilities $\frac14,\frac14,\frac12$
respectively; all random variables taking values in
$(-1,0)\bigcup(0,1)$ are completely random non-chain recurrent
variables; and the random variables taking values $1,-1,0$ with
positive probability and also taking values in $(-1,0)\bigcup(0,1)$
with positive probability are partly random chain recurrent
variables.
\end{example}

\begin{example}\rm
Let the probability space $(\Omega,\mathcal F,\mathbb P)$ be given
by $\Omega=S^1$ with $\mathcal F=\mathcal B(\Omega)$, and $\mathbb
P$ the Lebesgue measure. Put $\theta_t\omega=\omega+t$. Let the
state space $X=S^1$, too. Define a random homeomorphism
$\psi(\omega): S^1\rightarrow S^1$ by $\psi(\omega)x=x+\omega$. We
define an RDS by
\[
\varphi(t,\omega)=\psi(\theta_t\omega)\circ\varphi_0(t)\circ\psi^{-1}(\omega),
\]
where $\varphi_0$ is the flow on $S^1$ determined by the equation
\[
\dot{x}=-\cos x.
\]
Then the RDS $\varphi$ has no non-trivial attractor besides
$X,\emptyset$, for the details see \cite{Cra}. Hence by our main
theorem, the random chain recurrent set is $X$, i.e. all random
variables are random chain recurrent.
\end{example}

\section{Some discussions}

We know very well that there are several nonequivalent definitions
of random attractors for random dynamical systems, see \cite{Sc} for
instance. Pull-back attractors were introduced and studied by Crauel
and Flandoli \cite{Cra5}, Crauel, Debussche, and Flandoli
\cite{Cra3}, Schmalfuss \cite{Sch1,Sch2} and others. Ochs in
\cite{Och} firstly introduced  random weak attractors, where `weak'
means that the convergence to attractor is in probability instead of
usual almost sureness. Another kind of attractor is forward
attractor, which is in contrast to pull-back attractor and whose
convergence is almost sure convergence in contrast to weak
attractor's convergence in probability. As stated in \cite{Ar2}, the
choice of convergence in probability makes the forward and pull-back
attractors equivalent. So the authors adopted the weak attraction,
in fact they adopt the weak point attractor, in order to prove
Lyapunov's second method for RDS, for details see \cite{Ar2}. In
order to get the Morse theory for RDS, the authors in \cite{Cra4}
introduced the definition of attractor-repeller pair for RDS. They
defined an attractor to be the maximal invariant random compact set
inside its fundamental neighbourhood, i.e. a forward invariant
random open set, and defined the repeller corresponding to it to be
the complement of the basin of the attractor. For details, see
\cite{Cra4}.

In this paper, we adopt the definition \ref{UA}. Now we simply
discuss the relations between it and the old ones. It is obvious
that, for a random pre-attractor $U(\omega)$, the random attractor
$A(\omega)$ determined by it in definition \ref{UA} is a random
pull-back attractor in the universe
\[
\mathcal D=\{D(\omega)|~\overline{D(\omega)}\subset U(\omega)~{\rm
and}~D(\omega)~{\rm is ~measurable}\},
\]
therefore it is a weak set attractor, and therefore a weak point
attractor in $\mathcal D$. But it is obvious not necessarily a
forward attractor in $\mathcal D$. As mentioned above, the authors
in \cite{Cra4} introduced the definition of attractor-repeller pair
for RDS. It seems that their definition of attractor works for our
purpose, too. But when we prove lemma \ref{3.7}, we do not know how
to construct a fundamental neighbourhood $U(\omega)$ defined in
\cite{Cra4} to obtain the attractor. With respect to the relations
between our attractor and the attractor of \cite{Cra4}, the
difference is that in our definition, the random pre-attractor is
not necessarily forward invariant, while the fundamental
neighbourhood in \cite{Cra4} is forward invariant; our attractor is
pull-back attractor in the universe $\mathcal D$ defined above,
while the attractor in \cite{Cra4} is only a weak set attractor
inside its fundamental neighbourhood. Hence our attractor is not
weaker than the attractor in \cite{Cra4}, and vice versa.

In fact, besides these, our definition of attractor approximates
closely to Conley's deterministic definition of attractor. Conley in
\cite{Con} defined an attractor to be an invariant compact set which
is the omega-limit set of one of its neighbourhoods. A random
attractor defined by definition \ref{UA} is obvious the omega-limit
set ( a random set $A(\omega)$ is called the {\em omega-limit set}
of $D(\omega)$ if
$A(\omega)=\lim_{t\rightarrow\infty}\varphi(t,\theta_{-t}\omega)D(\theta_{-t}\omega)$
) of one of the pre-attractors determining it. Conversely, if an
invariant random compact set is the omega-limit set of one of its
neighbourhoods, then, by the definition of omega-limit set, we can
see that the neighbourhood must satisfies (\ref{at}), except that
the $T(\omega)$ in (\ref{at}) is not necessarily measurable, which
guarantees the measurability of $U(T(\omega))$, $A(\omega)$ etc.

In this paper we adopt the form
\[
X-\mathcal {CR}(\varphi)=\bigcup [B(A)-A]
\]
instead of the original form
\[
\mathcal {CR}(\varphi)=\bigcap [A\cup R]
\]
in our main theorem, where $R$ denotes the repeller corresponding to
the attractor $A$. On the one hand the form we employ here is also
adopted by other authors, see \cite{Hu0,Hu1,Hu2} for instances. In
fact we are inspired by these references. On the other hand, the
random pre-attractor is not necessarily forward invariant, which
restricts us to obtain the invariance of the basin of attractor as
in \cite{Cra4}. Hence we can not obtain the invariance of
$X-B(A)(\omega)$, but it is forward invariant, see lemma \ref{Ba}.
So we can not define $X-B(A)(\omega)$ to be the repeller
corresponding to the attractor $A(\omega)$ as in \cite{Cra4}, for a
repeller should be invariant.

With respect to the measurability, we find it not appropriate to
confine us to considering $\mathcal F$-measurability only. Since it
is easy to see that $U(T(\omega))$, $B(A)(\omega)$, $A(\omega)$,
$d(\omega)$ etc in this paper are not measurable with respect to
$\mathcal F$. Therefore to serve our purpose, we have to allow
measurability with respect to other $\sigma$-algebra, in fact we
allow $\mathcal F^u$-measurability. This treatment is also adopted
by others, see \cite{Cra1,Cra2,Cra3,Cra5,Chu} for instance.
Certainly we can choose ${\bar{\mathcal F}}^{\mathbb P}$, the
$\mathbb P$-completion $\sigma$-algebra of $\mathcal F$, in order
that the above variables are measurable. Of course, when $\mathbb
T=\mathbb Z$, i.e. the RDS in consideration is discrete, or
$\mathcal F$ is complete with respect to the probability measure
$\mathbb P$ on base space $(\Omega,\mathcal F)$, i.e. $\mathcal
F={\bar{\mathcal F}}^{\mathbb P}$, all the random variables are
measurable with respect to $\mathcal F$ as usual. Whence in this
case, only considering $\mathcal F$-measurability as usual is enough
to obtain our results.

\vskip 5mm \noindent{\bf\Large Acknowledgements}

The author expresses his sincere thanks to Professor Yong Li for his
instructions and many invaluable suggestions. The author is very
grateful to Professor Youqing Ji and Qingdao Huang for helpful
discussions. The author is also grateful to the anonymous referees
for their careful reading the manuscript and helpful suggestions,
one of whom has pointed out an important oversight and some
confusing notations.


\begin{thebibliography}{xx}

\bibitem{Ar1}
Arnold L 1998 \it Random Dynamical Systems \rm (Berlin Heidelberg
New York: Springer-Verlag)

\bibitem{Ar2}
Arnold L and Schmalfuss B 2001 Lyapunov¡¯s Second Method for Random
Dynamical Systems \it J. Diff. Equ. \bf 177 \rm 235-65

\bibitem{Cas}
Castaing C and Valadier M 1977 \it Convex Analysis and Measurable
Multifunctions \rm(\it Lec. Notes in Math. \rm vol 580) (Berlin:
Springer-Verlag)

\bibitem{Chu}
Chueshov I 2002 \it Monotone Random systems Theory and Applications
\rm(\it Lec. Notes in Math. \rm vol 1779) (Berlin Heidelberg:
Springer-Verlag)


\bibitem{Con}
Conley C 1978 \it Isolated Invariant Sets and the Morse Index \rm
(\it Conf. Board Math. Sci. \rm  vol 38) (Providence: Amer. Math.
Soc.)


\bibitem{Cra1}
Crauel H 1999 Global random attractors are uniquely determined by
attracting deterministic compact sets \it Ann. Mat. Pura Appl. (IV)
\bf 176 \rm  57-72

\bibitem{Cra2}
Crauel H 2001 Random Point Attractors versus Random Set Attractors
\it J. London Math. Society \bf 63 \rm (2) 413-27

\bibitem{Cra}
Crauel H 2002 A uniformly exponential random forward attractor which
is not a pullback attractor \it Arch. Math. \bf 78 \rm 329-36

\bibitem{Cra3}
Crauel H, Debussche A and Flandoli F 1997 Random Attractors \it J.
Dyn. Diff. Equ. \bf 9 \rm 307-41

\bibitem{Cra4}
Crauel H, Duc L H and Siegmund S 2004 Towards a Morse heory for
random dyanamical systems \it Stochastics and Dynamics \bf 4 \rm (3)
277-96

\bibitem{Cra5}
Crauel H and Flandoli F 1994 Attractors for random dynamical systems
\it Probab. Theory Related Fields \bf 100 \rm 365-93

\bibitem{Fra}
Franks J 1988 A variation of Poincar\'e-Birkhoff theorem {\it
Hamiltonian Dynamical Systems} (\it Contemporary Mathematics \rm vol
81) (Providence: Amer. Math. Soc.) pp. 111-17

\bibitem{Hu0}
Hurley M 1991 Chain recurrence and attraction in non-compact spaces
\it Ergod. Th. Dynam. Sys. \bf 11 \rm 709-29

\bibitem{Hu1}
Hurley M 1992 Noncompact chain recurrence and attraction \it Proc.
Amer. Math. Soc. \bf 115 \rm 1139-48

\bibitem{Hu2}
Hurley M 1995 Chain recurrence, semiflows, and gradients \it J. Dyn.
Diff. Equ. \bf 7 \rm 437-56

\bibitem{Klo}
Kloeden P E and Platen E 1992 \it Numerical Solution of Stochastic
Differential Equations \rm (Berlin Heidelberg New York:
Springer-Verlag)

\bibitem{Och}
Ochs G 1999 \it Weak Random Attractors \rm Institut f\"{u}r
Dynamische Systeme, Universit\"{a}t Bremen Report \bf 449 \rm

\bibitem{Sc}
Scheutzow M 2002 Comparison of various concepts of a random
attractor: A case study \it Arch. Math. \bf 78 \rm 233-40

\bibitem{Sch1}
Schmalfuss B 1992 Backward cocycles and attractors for stochastic
differential equations In: Reitmann V, Riedrich T and Koksch N
(Eds.) \it International Seminar on Applied Mathematics -- Nonlinear
Dynamics: Attractor Approximation and Global Behaviour \rm (Teubner,
Leipzig) pp. 185-92

\bibitem{Sch2}
Schmalfuss B 1997 The random attractor of the stochastic Lorenz
system \it ZAMP \bf 48 \rm 951-75



\end{thebibliography}
\end{document}